\documentclass[12pt,twoside]{extarticle}

\usepackage[utf8]{inputenc}
\usepackage{lipsum}    
\usepackage[utf8]{inputenc}
\usepackage[T1]{fontenc}
\usepackage[english]{babel}
\usepackage{geometry}
\usepackage{fancyhdr}
\usepackage{titlesec}
\usepackage{amssymb}\usepackage{amsthm}\usepackage{amsmath}\usepackage{mathrsfs}\usepackage{halloweenmath}
\usepackage[]{enumitem}
\usepackage{comment}
\usepackage[dvipsnames]{xcolor}
\usepackage{microtype}
\usepackage{hyperref}
\usepackage{mathabx}
\usepackage{mlmodern} 
\usepackage{moresize}
\usepackage{footmisc}
\usepackage{amsmath}
\usepackage{titling}
\usepackage[dvipsnames]{xcolor}

\setenumerate{itemsep=0pt,topsep=2pt,parsep=1pt}
\setitemize{itemsep = 0pt, topsep = 0pt, parsep = 0pt}
\hypersetup{colorlinks,linkcolor={red},citecolor={blue}}

\renewcommand{\thesection}{\arabic{section}}

\newcommand{\uc}{\MakeLowercase}
\newcommand{\siz}{\fontsize{13}{\baselineskip}\selectfont}
\newcommand{\sizz}{\fontsize{13.5}{\baselineskip}\selectfont}
\newcommand{\sizzz}{\fontsize{11}{\baselineskip}\selectfont}
\newcommand{\bibsiz}{\fontsize{9}{0}\selectfont}

\titleformat{\chapter}
    {\Large \centering \sc \vspace*{0cm}}
    {Chapter \thechapter:}
    {5pt}
    {}
    [\vspace{1cm} \titlerule]
\titleformat{\section}
    {\large \centering \bfseries} 
    {\thesection.} 
    {5pt} 
    {} 
\titleformat{\subsection}
    [runin]
    {\bfseries}
    {\thesubsection}
    {5pt}
    {\addperiod}
\newcommand{\addperiod}[1]{#1.}
\renewenvironment{proof}{{\itshape Proof.}}{\hfill\qedsymbol\vspace{8pt}}

\newtheoremstyle{defbfnote}
    {8pt} 
    {8pt} 
    {} 
    {\parindent} 
    {\scshape} 
    {.} 
    {.5 em} 
    {\thmname{#1}\thmnumber{ #2}\thmnote{ (\uc{\siz #3})}}
\newtheoremstyle{bfnote}
    {8pt} 
    {8pt} 
    {\itshape} 
    {\parindent} 
    {\scshape} 
    {.} 
    {.5 em} 
    {\thmname{#1}\thmnumber{ #2}\thmnote{ (\uc{\siz #3})}}

\theoremstyle{bfnote}
\newtheorem{theorem}{\uc{\siz Theorem}}
\newtheorem{claim}[theorem]{\uc{\siz Claim}}

\newtheorem{lemma}[theorem]{\uc{\siz Lemma}}

\theoremstyle{defbfnote}

\newtheorem{remark}[theorem]{\uc{\siz Remark}}

\newcommand{\Lcorner}{\mathrel{\makebox[7pt][c]{\rule{0.4pt}{5.75pt}\rule{4.5pt}{.4pt}}}}
\renewcommand{\llcorner}{{\Lcorner}}

\setlist[enumerate,1]{label={(\roman*)}}
\setlist[itemize,1]{label = {$-$}}
\DeclareMathAlphabet{\pzocal}{OMS}{zplm}{m}{n}
\DeclareMathAlphabet{\mathcal}{OMS}{cmsy}{m}{n}

\DeclareFontFamily{U}{ntxmia}{}
\DeclareFontShape{U}{ntxmia}{m}{it}{<-> ntxmia }{}
\DeclareFontShape{U}{ntxmia}{b}{it}{<-> ntxbmia }{}
\DeclareSymbolFont{lettersA}{U}{ntxmia}{m}{it}
\SetSymbolFont{lettersA}{bold}{U}{ntxmia}{b}{it}

\usepackage{amssymb,tikz}

\newcommand{\Xint}[1]{%
  \mathchoice
  {\XXint\displaystyle\textstyle{#1}}%
  {\XXint\textstyle\scriptstyle{#1}}%
  {\XXint\scriptstyle\scriptscriptstyle{#1}}%
  {\XXint\scriptscriptstyle\scriptscriptstyle{#1}}%
  \!\int
}
\newcommand{\XXint}[3]{%
  \setbox0=\hbox{$#1{#2#3}{\int}$}%
  \vcenter{\hbox{$#2#3$}}\kern-.5\wd0
}
\newcommand{\fint}{\Xint-}

\renewcommand{\mathbb}{\mathbf}

\renewcommand{\P}{\mathcal{P}}

\newcommand{\R}{\mathbb{R}}

\newcommand{\Rn}{{{\mathbb{R}}^n}}
\newcommand{\C}{{{\mathscr{C}}}}
\renewcommand{\L}{{\mathscr{L}}}

\newcommand{\spt}{\text{\normalfont spt}}

\newcommand{\tr}{\text{\normalfont tr}}
\newcommand{\LS}{\text{\normalfont LSF}}
\newcommand{\m}{\mathbf{m}}
\newcommand{\var}{\text{\normalfont cov}}

\newcommand{\Symp}{\text{\normalfont Sym}^+}
\renewcommand{\Bar}{\mathscr{B}}
\renewcommand{\b}{\mathbf{b}}
\renewcommand{\theta}{\vartheta}
\renewcommand{\epsilon}{\varepsilon}
\renewcommand{\phi}{\varphi}
\newcommand{\blambda}{{\boldsymbol{\lambda}}}
\renewcommand{\O}{\text{\normalfont O}}
\newcommand{\argmin}{\operatorname*{\text{\normalfont argmin}}}
\newcommand{\argmax}{\operatorname*{\text{\normalfont argmax}}}

\makeatletter
\def\thanks#1{\protected@xdef\@thanks{\@thanks
        \protect\footnotetext{#1}}}
\makeatother

\title{\uppercase{\bfseries \large On the closure under barycenters of location-scatter families with a normalized generator}}
\author{{\sc pietro aldrigo\textsuperscript{\textdagger}} \thanks{\textsuperscript{\textdagger}{\sc Universit\"at Bern, Mathematisches Institut}. Email address: \texttt{pietro.aldrigo@unibe.ch} 
The author is supported by the Swiss National Science Foundation, grant number 200021-228012. }}
\date{ }

\begin{document}
\maketitle
\vspace{-1.2cm}
\thispagestyle{empty} 

\begin{abstract}
We provide an example of a normalized absolutely continuous probability whose location-scatter family (LSF) does not contain all the 2-Wassertein geodesics passing through its elements. In particular, we find a counterexample for Corollary 4.5 in \cite{AlvarezEsteban2016}. Finally, we prove that LSFs with a (possibly singular) radially symmetric generator are closed under 2-Wasserstein barycenters. 
\end{abstract}

\section{Introduction}

The quadratic Wasserstein space $(\P_2(\Rn),W_2)$ is a geodesic metric space, and its geodesics are described by displacement interpolation along optimal transport plans \cite{AmbrosioGigliSavare2008,McCann1997,Villani2009}. A subset $\mathcal K\subseteq\P_2(\Rn)$ is said to be geodesically convex if for any two of its elements $\mu_1,\mu_2\in \mathcal{K}$ all the $W_2$-geodesic joining $\mu_1$ and $\mu_2$ are contained in $\mathcal K$. A natural multi-marginal generalization of this property is closure under Wasserstein barycenters. Given $\mu_1,\dots,\mu_N\in\P_2(\Rn)$ and weights $\blambda:=(\lambda_1,\ldots,\lambda_N)\in [0,1]^N$ with $\sum_{j=1}^N\lambda_j=1$, a Wasserstein barycenter is a minimizer
\begin{align*}
    \bar\mu\in\operatorname*{argmin}_{\nu\in\mathcal P_2(\mathbb R^n)}
    \sum_{i=1}^N\lambda_i W_2^2(\mu_i,\nu)=:\Bar_\blambda(\mu_1,\dots,\mu_N).
\end{align*}
Introduced by Agueh and Carlier \cite{AguehCarlier2011}, this construction recovers McCann's displacement interpolation when $N=2$, while for $N\geq3$ it is naturally related to the multi-marginal optimal transport problem with barycentric quadratic cost (see Theorem \ref{thAG} below) studied by Gangbo and {\'S}wi{\k{e}}ch \cite{GangboSwiech1998},  see also \cite{Gelbrich1990,Pass2015}.

A family $\mathcal{F}\subseteq \P_2(\Rn)$ is closed under barycenters if $\Bar_\blambda(\mu_1,\dots,\mu_N)\subseteq \mathcal{F}$ for every $\mu_1,\dots,\mu_N\in \mathcal{F}$ and weights $\blambda$. While it is clear that closure under barycenters implies geodesic convexity, the converse implication does not hold in general (see e.g. \cite{EmamiPass2025}). 

An important class of families for which closure under barycenters has been studied explicitly is that of location-scatter families. Given $\mu_0\in\mathcal P_{2}(\Rn)$ one defines the location-scatter family generated by $\mu_0$ as
\begin{align*}
 \LS(\mu_0):= \left\{(f_{A,m})_\#\mu_0: f_{A,m}(x):=Ax + m,\,A\in \Symp(\Rn),\,m\in \Rn\right\},
\end{align*}
where $\Symp(\Rn)\subseteq \Rn\otimes \Rn$ denotes the family of positive symmetric matrices. Location-scatter families are not in general closed under barycenters. However, \emph{under suitable assumptions on the generator $\mu_0$}, for any choice of weights $\blambda$, the Wasserstein barycenter of measures $\{(f_{A_j,m_j})_\#\mu_0\}_{1\leq j\leq N}$ is unique and admits an explicit expression in terms of the push-forward of $\mu_0$ through an affine map $f_{A,m}(x):= Ax 
+m$, where $A\in \Symp(\Rn)$ depends only on $\blambda$ and the matrices $A_1,\dots,A_N$, and $m = \sum_{j=1}^N\lambda_j m_j$; in particular, in this setting, $\LS(\mu_0)$ is closed under under barycenters. This is for instance the case for Gaussians distributions \cite{AguehCarlier2011,AlvarezEsteban2016}.
 
Recall that for a probability measure $\mu$ on $\Rn$ with finite second moments, the \emph{mean of $\mu$} and \emph{covariance matrix of $\mu$} are respectively defined as 
\begin{align*}
\m(\mu):=\int_\Rn x\,d\mu(x)\in \Rn\,\, \text{and}\,\, \var(\mu):=\int_{\Rn} (x-\m(\mu))\otimes (x-\m(\mu))\,d\mu(x)\in \Rn\otimes \Rn.
\end{align*}
In Corollary 4.5 of \cite{AlvarezEsteban2016}, the authors stated the following result.

\begin{claim}\label{th}
	Let $\mu_0\in \P_2(\Rn)$ be an absolutely continuous probability measure, with zero mean and identity covariance. Then the family $\LS(\mu_0)$ is closed under Wasserstein barycenters. 
\end{claim}

The aforementioned result has also been cited in several other works, a non-exhaustive list being \cite{AlvarezEsteban2018,del2019robust,dusson2026wasserstein,korotin2022wasserstein,korotin2021continuous}.
In this paper however, we exhibit a counterexample to this statement. More precisely, we prove the existence of a normalized absolutely continuous probability measure $\nu_0\in \P_2(\R^2)$ such that there exist two elements $\nu_1,\nu_2\in \LS(\nu_0)$ with the property $\Bar_{(1/2,1/2)}(\nu_1,\nu_2)\cap \LS(\nu_0)=\varnothing$; hence, proving that $\LS(\nu_0)$ is not even geodesically convex.

The result stated above, while false, turns out to be not “far from being true”. Indeed we have the following.

\begin{theorem}\label{th_true}
	Let $\mu_0\in \P_2(\Rn)$ be a radially symmetric probability measure. Then $\LS(\mu_0)$ is closed under Wasserstein barycenters. More precisely, for any $N\geq 2$, $\blambda = (\lambda_1,\dots,\lambda_N)\in [0,1]^N$ with $\sum_{j=1}^N \lambda_j=1$ and $\mu_j = (f_{A_j,m_j})_\#\mu_0\in \LS(\mu_0)$, for $1\leq j\leq N$, then $\Bar_\blambda(\mu_1,\dots,\mu_N)$ contains one and only one element, and this is $\mu_{\blambda}:= (f_{X_\blambda,m_\blambda})_\#\mu_0$, where $m_\blambda :=\sum_{j=1}^N\lambda_j m_j$ and $X_\blambda$ is the unique symmetric positive solution to  
	\begin{align}\label{def_X}
	X_\blambda^2 = \sum_{j=1}^N \lambda_j (X_\blambda A_j^2 X_\blambda)^{\frac{1}{2}}.
	\end{align}
\end{theorem}

\textbf{Structure of the paper.} In Section \ref{sec_Notation} we establish the notation used throughout this work. In Section \ref{sec_Counter} we provide a counterexample to Claim \ref{th}. Finally, in Section \ref{sec_Proof} we prove Theorem \ref{th_true}.

\medskip

\textbf{Acknowledgments.} 
I would like to thank my supervisor, Prof. Dr. Zolt\'{a}n Balogh, for the many fruitful discussions we have had during the development of this work.

\section{Notation and preliminaries}\label{sec_Notation}

Let us denoted by $\P_2(\Rn)$ the family of Borel probability measures on $\Rn$ with finite second moments. 
A Borel probability $\mu\in \P_2(\Rn)$ is said to be \emph{normalized} if $\m(\mu)=0\in \Rn$ and $\var(\mu)=\mathbf{1}_n$.

We denote by $\O(\Rn)$ the family of matrices $R\in \Rn\otimes \Rn$ such that $RR^t = \mathbf{1}_n = R^tR$, and we call its elements \emph{orthogonal matrices}. We say that a probability measure $\mu_0\in \P_2(\Rn)$ is \emph{radially symmetric} if $\mu_0 = R_\#\mu_0$ for every every $R\in \O(\Rn)$, where the elements of $\Rn\otimes \Rn$  are identified with linear functions $\Rn\to \Rn$. For each $m\in \Rn$, we write $\tau_m$ to denote the translation $x\mapsto x+m$. If $f_1,\dots,f_N:\Rn\to \Rn$, we define $\bigtimes_{j=1}^N f_j:\Rn\to (\Rn)^N$ and $\bigotimes_{j=1}^N f_j:(\Rn)^N\to (\Rn)^N$ as
\begin{align*}
\bigtimes_{j=1}^Nf_j(x):= (f_1(x),\dots,f_N(x))\quad \text{and}\quad \bigotimes_{j=1}^Nf_j(x_1,\dots,x_N):= (f_1(x_1),\dots,f_N(x_N))
\end{align*}
for all $x,x_1,\dots,x_N\in \Rn$.

For any $N\geq 2$, $\Sigma^{N-1}\subseteq \R^N$ is the set of elements $\mathbf{\blambda} = (\lambda_1,\dots,\lambda_N)\in [0,1]^N$ such that $\sum_{j=1}^N\lambda_j = 1$,. The elements of $\Sigma^{N-1}$ will be called \emph{weights}. 
Given weights $\blambda \in \Sigma^{N-1}$, we define $\b_{\mathbf\blambda} : (\Rn)^N\to \Rn$ as
\begin{align*}
\b_\blambda(x_1,\dots,x_N) := \sum_{j=1}^N \lambda_j x_j.
\end{align*}
If the weights $\mathbf\blambda$ satisfy $\lambda_1=\dots =\lambda_N = \lambda = 1/N$, we may write $\Bar_{\lambda}$ and $\b_\lambda$ instead of $\Bar_{\mathbf\blambda}$ and $\b_{\mathbf\blambda}$. Finally, a curve $\mu_{(\cdot)}:[0,1]\to \P_2(\Rn)$ is a $W_2$-geodesic if $\mu_{(\cdot)}$ is an isometric embedding of $[0,1]$ into $\P_2(\Rn)$.

For any $1\leq j\leq N$, we define $\pi_j:(\Rn)^N\to \Rn$ to be the   projection $\pi_j(x_1,\dots,x_N):= x_j$. Given $\mu_1,\dots,\mu_N\in \P_2(\Rn)$, we say that $\gamma\in \P_2((\Rn)^N)$ is a \emph{coupling of $\mu_1,\dots,\mu_N$}, and write $\gamma \in \Pi(\mu_1,\dots,\mu_N)$ if $(\pi_j)_\#\gamma = \mu_j$ for all $1\leq j\leq N$.

We recall the characterization of barycenters in terms of optimal plans for the multi-marginal transport problem due to Agueh and Carlier \cite{AguehCarlier2011}.

\begin{theorem}[Agueh--Carlier]\label{thAG}
	Let $N\geq 2$, $\blambda\in \Sigma^{N-1}$ and $\mu_1,\dots,\mu_N\in \P_2(\Rn)$. Then
	\begin{align}\label{thAG_eq}
	\min_{\mu\in \P_2(\Rn)}\sum_{j=1}^N \lambda_j W_2^2(\mu_j,\mu) = \min_{\gamma\in \Pi(\mu_1,\dots,\mu_N)}\int_{(\Rn)^N}\sum_{1\leq i,j\leq N} \lambda_i\lambda_j \frac{|x_i-x_j|^2}{2}\,d\gamma(x_1,\dots,x_N).
	\end{align}
	Moreover, $\mu_*\in \P_2(\Rn)$ realizes the minimum on the left-hand side of \eqref{thAG_eq} if and only if there exists $\gamma_*\in \Pi(\mu_1,\dots,\mu_N)$ minimizer of the right-hand side of \eqref{thAG_eq} such that $\mu_* = (\b_\blambda)_\#\gamma_*$. 
\end{theorem}

Choosing $N=2$ in Theorem \ref{thAG} yields the following consequences.
\begin{enumerate}
\item  $\mu_* \in \Bar_{(1-t,t)}(\mu_0,\mu_1)$ for some $0\leq t\leq 1$ if and only if there exists a $W_2$-geodesics ${\mu}_{(\cdot)}$ that joins $\mu_0$ and $\mu_1$, and $\mu_* = {\mu}_t$; and in particular, closure under barycenters is a stronger property than geodesic convexity.
\item If $\mu$ is absolutely continuous, Brenier's theorem guarantees existence and uniqueness of Wasserstein barycenter $\Bar_{(1-\lambda,\lambda)}(\mu,\nu)$ for any $\nu\in \P_2(\Rn)$ and any weights $(1-\lambda,\lambda)$.
\end{enumerate}

\section{A counterexample to Claim \ref{th}}\label{sec_Counter}

Consider the probability measure $\mu_0\in \P_2(\R^2)$ defined by
\begin{align}\label{def_mu0}
\mu_0 := \frac{6}{10}\delta_{p_1}+ \frac{3}{10}\delta_{p_2}+ \frac{1}{10}\delta_{p_3},\quad  p_1:=\begin{pmatrix}
\sqrt{5}/3 \\ 1/3
\end{pmatrix},
\quad p_2 := \begin{pmatrix}
-2\sqrt{5}/3 \\ 1/3
\end{pmatrix},
\quad p_3 := \begin{pmatrix}
0 \\ - 3 
\end{pmatrix}.
\end{align}
Direct computations show that $\m(\mu_0) = 0$ and $\var(\mu_0) = \mathbb{1}_2$. Therefore, $\mu_0$ is a normalized probability measure. Consider the two symmetric matrices 
\begin{align}\label{def_matrices}
L:= \begin{pmatrix}
1 & 3 \\ 3 & 10
\end{pmatrix},
\quad M:= \begin{pmatrix}
1 & -3 \\ -3 & 10
\end{pmatrix}.
\end{align}
Then $\tr\, L = \tr\, M >0$ and $\det \, L = \det\, M = 1$. Therefore $L,M \in \Symp(\Rn)$ and $\mu_L := L_\#\mu_0,\mu_M := M_\#\mu_0\in \LS(\mu_0)$. Define $q_j := Lp_j$ and $r_j:= Mp_j$ for all $1\leq j\leq 3$. Then 
\begin{align*}
\mu_L = \frac{6}{10}\delta_{q_1}+ \frac{3}{10}\delta_{q_2}+ \frac{1}{10}\delta_{q_3} \quad \text{and}\quad \mu_M = \frac{6}{10}\delta_{r_1}+ \frac{3}{10}\delta_{r_2}+ \frac{1}{10}\delta_{r_3},
\end{align*}
and we have
\begin{align}\label{discrete_1}
|q_1 - r_2 | < |q_1-r_1|< |q_1-r_3|, \quad |q_2 - r_1 | < \min_{j\neq 1}|q_2 - r_j|, \quad |q_3-r_3| < \min_{j\neq 3}{|q_3-r_j|}.
\end{align}
From \eqref{discrete_1} we infer that the unique element of $\Pi(\mu_L,\mu_M)$ with cyclically monotone support is
\begin{align}\label{def_gamma*}
\gamma_* := \frac{3}{10}\delta_{(q_1,r_1)} + \frac{3}{10}\delta_{(q_1,r_2)} + \frac{3}{10}\delta_{(q_2,r_1)} + \frac{1}{10}\delta_{(q_3,r_3)}.
\end{align} 
Hence, $\gamma_*$ is the unique $W_2$-optimal transport plan from $\mu_L$ to $\mu_M$. Recalling Theorem \ref{thAG}, we deduce that the unique element of $\Bar_{1/2}(\mu_L,\mu_M)$ is 
\begin{align*}
\mu_* := (\b_{1/2})_\#\gamma_* = \frac{3}{10}\delta_{\frac{q_1+r_1}{2}} + \frac{3}{10}\delta_{\frac{q_1+r_2}{2}} + \frac{3}{10}\delta_{\frac{q_2+r_1}{2}} + \frac{1}{10}\delta_{\frac{q_3+r_3}{2}},
\end{align*}
and therefore the support $\mu_*$ consists of $4$ distinct points of the plane.

Fix $\epsilon>0$ small and consider the points $p_j^\epsilon:= \lambda_\epsilon p_j$, $1\leq j\leq 3$, where $\lambda_\epsilon:=\sqrt{1-\epsilon^2/4}$. Let $B_\epsilon^j$ be the ball with radius $\epsilon$ and center $p_j^\epsilon$ for every $1\leq j\leq 3$ and define
\begin{align*}
\mu_0^\epsilon := \frac{6}{10}\frac{\L^2\llcorner B_\epsilon^1}{\pi \epsilon^2} + \frac{3}{10}\frac{\L^2\llcorner B_\epsilon^2}{\pi \epsilon^2} + \frac{1}{10}\frac{\L^2\llcorner B_\epsilon^3}{\pi \epsilon^2},
\end{align*}
where, for each Borel set $\Omega\subseteq \R^2$, $\L^2\llcorner \Omega$ denotes the restriction of the 2-dimensional Lebesgue measure to $\Omega$.
Then it is clear that $\m(\mu_0^\epsilon) = \lambda_\epsilon \m(\mu_0) = 0$,
and, recalling that $\var(\mu_0)=\mathbf{1}_2$, and that for any $z=(z_1,z_2)\in \R^2$ and any $r>0$ we have
\begin{align*}
\fint_{B_r(z)} y_1^2\,dy = z_1^2 + \frac{r^2}{4}, \quad \fint_{B_r(z)} y_2^2\,dy = z_2^2 + \frac{r^2}{4}, \quad \fint_{B_r(z)} y_1y_2\,dy = z_1z_2,
\end{align*}
we also deduce that $\var(\mu_0^\epsilon) = \mathbb{1}_2$.
Hence, $\mu_0^\epsilon$ is a normalized absolutely continuous probability measure as long as $\epsilon>0$ is small enough. Observe also that the map $T^\epsilon:\R^2\to \R^2$ defined $\mu_0^\epsilon$-almost everywhere by the condition
$T^\epsilon(B_\epsilon^j):=p_j\quad \forall 1\leq j\leq 3,$
is a transport map from $\mu_0^\epsilon$ to $\mu_0$, and 
\begin{align*}
\lim_{\epsilon\downarrow 0}\int_{\R^2} |T^\epsilon(x)-x|^2\,d\mu_0^\epsilon(x) = 0.
\end{align*}
Therefore $\mu_0^\epsilon\to \mu_0$ in $\P_2(\Rn)$ as $\epsilon\downarrow 0$.

Define  the push-forwards $\mu_L^\epsilon:= L_\#\mu_0^\epsilon$ and $\mu_M^\epsilon:=M_\#\mu_0^\epsilon$ for each $\epsilon>0$ sufficiently small. The measures $\mu_L^\epsilon$ and $\mu_M^\epsilon$ are compactly supported and absolutely continuous; therefore, by Brenier's theorem and Theorem \ref{thAG}, there exists one unique $W_2$-optimal transport plan $\gamma^\epsilon_*\in \Pi(\mu_L^\epsilon,\mu_M^\epsilon)$ and one unique element $\mu_*^\epsilon\in \Bar_{1/2}(\mu_L^\epsilon,\mu^\epsilon_M)$. As $\mu_L^\epsilon\xrightarrow[]{\text{narrow}} \mu_L$ and $\mu_M^\epsilon\xrightarrow[]{\text{narrow}} \mu_M$ as $\epsilon\downarrow 0$, by uniqueness of the optimal coupling between $\mu_L$ and $\mu_M$, we deduce that $\gamma_*^\epsilon$ converges narrowly to the measure $\gamma_*$ defined in \eqref{def_gamma*}. By continuity of the map $\b_{1/2}:\R^2\times \R^2\to \R^2$, we also have
\begin{align}\label{eq_narrow}
\mu_*^\epsilon\xrightarrow[]{\text{narrow}} \mu_* = \frac{3}{10}\delta_{\frac{q_1+r_1}{2}} + \frac{3}{10}\delta_{\frac{q_1+r_2}{2}} + \frac{3}{10}\delta_{\frac{q_2+r_1}{2}} + \frac{1}{10}\delta_{\frac{q_3+r_3}{2}}.
\end{align}

Assume $\LS(\mu_0^\epsilon)$ to be geodesically convex for every $\epsilon$ in some right neighborhood of $0$. Then there exist families $\{A_\epsilon\}_{\epsilon>0}\subseteq \Symp(\R^2)$ and $\{m_\epsilon\}_{\epsilon>0}\subseteq \R^2$ such that $\mu_*^\epsilon = (f_\epsilon)_\#\mu_0^\epsilon$, where $f_\epsilon(x):= A_\epsilon x + m_\epsilon$, for all $\epsilon>0$ small enough. Since $\mu_0^\epsilon$ is normalized, we infer
\begin{align}\label{eq_moments}
m_\epsilon = \m(\mu_*^\epsilon) = \int_{\R^2}x\,d\mu_*^\epsilon(x) \quad \text{and}\quad A_\epsilon^2 = \var(\mu_*^\epsilon) = \int_{\R^2}x\otimes x\,d\mu_*^\epsilon(x)- m_\epsilon\otimes m_\epsilon.
\end{align}
Since the family $\{\mu_*^\epsilon\}_{\epsilon>0}$ is equi-compactly supported, narrow convergence of $\mu_*^\epsilon$ to $\mu_*$ implies the existence of the limits as $\epsilon\downarrow 0$ of the integrals appearing in \eqref{eq_moments}; hence we can define
\begin{align*}
m:= \lim_{\epsilon\downarrow 0} \int_{\R^2}x\,d\mu_*^\epsilon(x) \quad \text{and}\quad  A:= \left(\lim_{\epsilon\downarrow 0}\int_{\R^2}x\otimes x\,d\mu_*^\epsilon(x)-m\otimes m\right)^{\frac{1}{2}},
\end{align*}
and $f(x):= Ax + m$. By continuity of the operator $(\cdot)^\frac{1}{2}$ on the space of non-negative symmetric matrices, we deduce that $f_\epsilon$ converges uniformly to $f$ on compact sets. Therefore, for any test function $\phi\in \C^0_b(\R^2)$, we have
\begin{align}\label{eq_convergence}
\begin{split}
\lim_{\epsilon\downarrow 0 }\int_{\R^2} \phi(x)\,d\mu_*^\epsilon(x) 
	&= \lim_{\epsilon\downarrow 0}\int_{\R^2} (\phi(f_\epsilon(x))-\phi(f(x)))\,d\mu_0^\epsilon(x) + \int_{\R^2}\phi(x)\,d(f_\#\mu_0)(x)\\
	&=\int_{\R^2}\phi(x)\,d(f_\#\mu_0)(x),
\end{split}
\end{align}
where the last equality of \eqref{eq_convergence} follows from equi-compactness of the supports of $\{\mu_0^\epsilon\}_{\epsilon>0}$ and Heine--Cantor's theorem. Combining \eqref{eq_narrow} and \eqref{eq_convergence}, we deduce that $\mu_* = f_\#\mu_0$, which is not possible because the support of $\mu_*$ consists of 4 distinct points. 

The above argument proves the existence of $\overline{\epsilon}>0$ such that $\LS(\mu_0^{\overline{\epsilon}})$ is not geodesically convex; hence $\LS(\mu_0^{\overline{\epsilon}})$ is a fortiori not closed under barycenters. 

\section{Proof of Theorem \ref{th_true}}\label{sec_Proof}

\newcommand{\bl}{\blambda}

\begin{lemma}\label{Lem_reduction}
	Let $N\geq 2$, $\blambda \in \Sigma^{N-1}$, $m_1,\dots,m_N\in \Rn$, $\mu_1,\dots,\mu_N\in \P_2(\Rn)$ and define $\nu_j := (\tau_{m_j})_\#\mu_j$ for all $1\leq j\leq N$. Then $\nu_\blambda \in \Bar_\blambda(\nu_1,\dots,\nu_N)$ if and only if there exists $\mu_\blambda\in \Bar_\blambda(\mu_1,\dots,\mu_N)$ such that $\nu_\blambda = (\tau_{m_\blambda})_\#\mu_\blambda$, where $m_\blambda := \sum_{j=1}^N \lambda_j m_j$.
\end{lemma}

\begin{proof}
	By virtue of Theorem \ref{thAG}, $\mu_\blambda \in \Bar_\blambda(\mu_1,\dots,\mu_N)$ if and only if there exists $\gamma_\blambda\in \Pi(\mu_1,\dots,\mu_N)$ that maximizes the function
	\begin{align}\label{eq_function}
	\gamma \mapsto \int_{(\Rn)^N}\sum_{1\leq i,j\leq N}\lambda_i\lambda_j\langle x_i,x_j\rangle\,d\gamma(x_1,\dots,x_N)
	\end{align}
	in $\Pi(\mu_1,\dots,\mu_N)$ such that $\mu_\blambda = (\b_\blambda)_\#\gamma_\blambda$. 
	
	The map $\gamma \mapsto (\bigotimes_{j=1}^N\tau_{m_j})_\#\gamma$ is a bijection between $\Pi(\mu_1,\dots,\mu_N)$ and $\Pi(\nu_1,\dots,\nu_N)$, and $\gamma_\blambda$ maximizes the function in \eqref{eq_function} if and only if $\delta_\blambda:= (\bigotimes_{j=1}^N\tau_{m_j})_\#\gamma_\blambda$ maximizes the function 
	\begin{align*}
	\delta \mapsto \int_{(\Rn)^N}\sum_{1\leq i,j\leq N}\lambda_i\lambda_j\langle y_i,y_j\rangle\,d\delta(y_1,\dots,y_N)
	\end{align*}
	in $\Pi(\nu_1,\dots,\nu_N)$. This, together with Theorem \ref{thAG} and the observation that $\tau_{m_\bl}\circ \b_\bl = \b_\bl\circ \left(\bigotimes_{j=1}^N \tau_{m_j}\right)$, concludes the proof.

\end{proof}

	Let us finally remark that Theorem 6.1 of \cite{AguehCarlier2011} guarantees  existence and uniqueness of positive definite solutions for the matrix-equation  
	\begin{align*}
	X^2 = \sum_{j=1}^N \lambda_j (XA_j^2X)^{\frac{1}{2}},
	\end{align*}
	for any $A_1,\dots,A_N\in \Symp(\Rn)$ and $\bl\in \Sigma^{N-1}$. 
	
	\medskip
	{\itshape Proof of Theorem \ref{th_true}.} Thanks to Lemma \ref{Lem_reduction}, it is enough to prove the statement for $m_1 = \cdots = m_N = 0\in \Rn$. Also, by virtue of the above remark, the matrix $X_\bl\in \Symp(\Rn)$ is well-defined by \eqref{def_X}. Denote by $C_1,\dots,C_N$ and  $R_1,\dots,R_N$  the elements of $\Rn\otimes \Rn$ defined by
	\begin{align*}
	C_j:= X_\bl^{-1}(X_\bl A_j^2 X_\bl)^\frac{1}{2}X_{\bl}^{-1},\quad R_j := A_{j}^{-1}X_\blambda^{-1}(X_\bl A_j^2 X_\bl)^{\frac{1}{2}}\quad \forall 1\leq j\leq N.
	\end{align*}
	Observe that $C_1,\dots,C_N\in \Symp(\Rn)$ and, thanks to \eqref{def_X}, they satisfy the property 
	\begin{align}\label{eq_Cj}
	\sum_{j=1}^N \lambda_j C_j = X_\bl^{-1}\left(\sum_{j=1}^N \lambda_j (X_\bl A_j^2 X_\bl)^\frac{1}{2}\right)X_\bl^{-1}= \mathbf{1}_{n}.
	\end{align}
	Also, notice that 
	\begin{align*}
	R_jR^t_j = A_{j}^{-1}X_\blambda^{-1}(X_\bl A_j^2 X_\bl) X_\bl^{-1} A_j^{-1} = \mathbf{1}_n \quad \forall 1\leq j\leq N,
	\end{align*}
	which implies $R_1,\dots,R_N\in \O(\Rn)$, and
	\begin{align}\label{ajcjrjx}
	A_jR_j = C_jX_\lambda \quad \forall 1\leq j\leq N.
	\end{align}
	
	For any $x_1,\dots,x_N\in \Rn$ we have
	\begin{align*}
	\begin{split}
	\sum_{1\leq i,j\leq N}\lambda_i\lambda_j |x_i-x_j|^2 
	&= 2\sum_{j=1}^N \lambda_j|x_j|^2 -2\langle \sum_{i=1}^N\lambda_i x_i,\sum_{j=1}^N\lambda_jx_j\rangle\\
	&= 2\sum_{j=1}^N \lambda_j|x_j|^2 - 2 |\b_\bl(x_1,\dots,x_N)|^2.
	\end{split}
	\end{align*}
	Therefore,
	\begin{align}\label{argmin_argmax}
	\argmin_{\gamma\in \Pi(\mu_1,\dots,\mu_N)}\int_{(\Rn)^N} \sum_{1\leq i,j\leq N}\lambda_i\lambda_j |x_i-x_j|^2\,d\gamma(x_1,\dots,x_N) = \argmax_{\gamma\in \Pi(\mu_1,\dots,\mu_N)}\int_{(\Rn)^N} |\b_\bl|^2\,d\gamma.
	\end{align}
	Since $C_1,\dots,C_N$ are positive definite and \eqref{eq_Cj} holds, for every $x_1,\dots,x_N\in \Rn$, writing $b:=\b_\bl(x_1,\dots,x_N)$, we deduce
	\begin{align*}
	0\leq \sum_{j=1}^N \lambda_j \langle C_j^{-1}(x_j - C_jb),x_j - C_jb\rangle = \sum_{j=1}^N \lambda_j\langle C_j^{-1}x_j,x_j\rangle - |b|^2.
	\end{align*}
	This proves that
	\begin{align}\label{ineq_b}
	|\b_\bl(x_1,\dots,x_N)|^2 \leq \sum_{j=1}^N\lambda_j \langle C_j^{-1}x_j,x_j\rangle \quad \forall x_1,\dots,x_N\in \Rn,
	\end{align}
	and equality holds in \eqref{ineq_b} if and only if 
	\begin{align*}
	x_j = C_j\b_\bl(x_1,\dots,x_N)\quad \forall 1\leq j\leq N.
	\end{align*}
	
	Let $\gamma\in \Pi(\mu_1,\dots,\mu_N)$ be any arbitrary coupling. Integrating \eqref{ineq_b} with respect to $\gamma$ yields
	\begin{align}\label{eq_plan}
	\int_{(\Rn)^N}|\b_\bl|^2\,d\gamma \leq \sum_{j=1}^N\lambda_j \int_\Rn \langle C_j^{-1}x_j,x_j\rangle \,d\mu_j(x_j).
	\end{align}
	On the other hand, setting $\gamma_\bl := (\bigtimes_{j=1}^N A_jR_j)_\#\mu_0$, recalling that $\mu_0$ is rotationally symmetric and taking into account \eqref{def_X}, \eqref{eq_Cj} and \eqref{ajcjrjx}, we infer that
	\begin{align}\label{eq_optimalplan}
	\begin{split}
	\int_{(\Rn)^N}&|\b_\bl|^2\,d\gamma_\bl 
	= \int_{\Rn} \left|\sum_{j=1}^N\lambda_jC_jX_\bl x\right|^2\,d\mu_0(x) = \int_\Rn \langle X_\bl^2 x,x\rangle\,d\mu_0(x) \\
	&= \sum_{j=1}^N \lambda_j\int_\Rn \langle (X_\bl A_j^2X_\bl)^{\frac{1}{2}}x,x\rangle\,d\mu_0(x)
	= \sum_{j=1}^N \lambda_j\int_\Rn \langle X_\bl C_j X_\bl x,x\rangle\,d\mu_0(x)\\
	&= \sum_{j=1}^N \lambda_j\int_\Rn \langle C_j^{-1} C_j X_\bl x,C_j X_\bl x\rangle\,d\mu_0(x)
	= \sum_{j=1}^N \lambda_j\int_\Rn \langle C_j^{-1}A_jR_j x,A_j R_j x\rangle\,d\mu_0(x)\\ 
	&= \sum_{j=1}^N \lambda_j\int_\Rn \langle C_j^{-1} A_j x,A_j x\rangle\,d\mu_0(x)
	 = \sum_{j=1}^N \lambda_j\int_\Rn \langle C_j^{-1} x, x\rangle\,d\mu_j(x).
	\end{split} 
	\end{align}
	From \eqref{eq_plan} and \eqref{eq_optimalplan}, it follows that $\gamma_\bl$ belongs to $\argmax_{\gamma\in\Pi(\mu_1,\dots,\mu_N)}\int_{(\Rn)^N}|\b_\bl|^2\,d\gamma$. Moreover, by the characterization of the equality case in \eqref{ineq_b}, any other maximizer $\gamma_*\in \Pi(\mu_1,\dots,\mu_N)$ satisfies the property 
	\begin{align}\label{eq_*}
	x_j = C_j\b_\bl(x_1,\dots,x_N)\quad \forall 1\leq j\leq N\quad \gamma_*\text{-almost everywhere}.
	\end{align}
	
	Let $\gamma_*\in \Pi(\mu_1,\dots,\mu_N)$ be such that \eqref{eq_*} holds true, and let $\phi:\Rn\to [0,\infty)$ be any Borel-measurable function. Then  
	\begin{align}\label{eq_p}
	\begin{split}
	\int_\Rn \phi\,d(\b_\bl)_\#\gamma_* &= \int_{(\Rn)^N} \phi(\b_\bl(x_1,\dots,x_N))\,d\gamma_*(x_1,\dots,x_N)\\
	& = \int_{(\Rn)^N}\phi(C_1^{-1}x_1)\,d\gamma_*(x_1,\dots,x_N)\\
	&= \int_{\Rn}\phi(C_1^{-1}x_1)\,d\mu_1(x_1)
	 = \int_\Rn \phi(C_1^{-1}A_1x)\,d\mu_0(x)\\
	&= \int_\Rn \phi(C_1^{-1}A_1R_1x)\,d\mu_0(x)
	= \int_\Rn \phi\,d(X_\bl)_\#\mu_0.
	\end{split}
	\end{align}
	By Theorem \ref{thAG} it follows that every element of $\Bar_\bl(\mu_1,\dots,\mu_N)$ is of the from $(\b_\bl)_\#\gamma_*$, for some $\gamma_*\in \argmax_{\gamma\in \Pi(\mu_1,\dots,\mu_N)}\int_{(\Rn)^N}|\b_\bl|^2\,d\gamma$. Hence, we deduce from \eqref{eq_p} that $\Bar_\bl(\mu_1,\dots,\mu_N)$ contains one unique element, and this is $(X_\bl)_\#\mu_0 = \mu_\bl$, which belongs to $\LS(\mu_0)$.
	
	\hfill \qedsymbol
	
	\begin{remark}
	Observe that the argument provided in the previous proof can be used to show   that
	\begin{align}\label{eq_q}
	\argmin_{\gamma\in \Pi(\mu_1,\dots,\mu_N)}\int_{(\Rn)^N}\sum_{1\leq i,j\leq N}\lambda_i\lambda_j|x_i-x_j|^2\,d\gamma(x_1,\dots,x_N) = \left\{(\bigtimes_{j=1}^N A_jR_j)_\#\mu_0\right\}.
	\end{align}
	Indeed, we have already established that any element $\gamma_*$ belonging to the left hand-side must satisfy \eqref{eq_*}, and therefore 
	\begin{align}\label{gammastar}
	\gamma_* = (\bigtimes_{j=1}^N C_j)_\# (\b_\bl)_\#\gamma_*.
	\end{align}
	But \eqref{eq_p} states that any such measure must also satisfy 
	\begin{align*}
	(\b_\bl)_\# \gamma_* = (X_\bl)_\#\mu_0 = \mu_\bl.
	\end{align*}
	Hence, by \eqref{gammastar}, all the elements of the left hand-side of \eqref{eq_q} are of the form
	\begin{align*}
	\gamma_* = (\bigtimes_{j=1}^N C_j)_\#(X_\bl)_\#\mu_0 = (\bigtimes C_jX_\blambda)_\#\mu_0 = (\bigtimes_{j=1}^N A_jR_j)_\#\mu_0.
	\end{align*}
	\end{remark}

{\bibsiz
\bibliographystyle{acm}

}

\end{document}